# Energy-based dual-phase dynamics identification of clearance nonlinearities

Cristian López[1], Keegan J. Moore[1,2]*

[1]*Department of Mechanical and Materials Engineering, University of Nebraska-Lincoln, Lincoln, NE 68588*
[2]*Daniel Guggenheim School of Aerospace Engineering, Georgia Institute of Technology, Atlanta, GA 30332*

*Corresponding author: (K.J. Moore)
E-mail address: kmoore@gatech.edu

## ABSTRACT

The energy-based dual-phase dynamics identification (EDDI) method is a new data-driven technique for the discovery of equations of motion (EOMs) of strongly nonlinear single-degree-of-freedom (SDOF) oscillators. This research uses the EDDI method to obtain mathematical models for SDOF systems with clearance nonlinearities. The first key aspect of the EDDI method is that it relates the kinetic energy of the system to the dissipated energy and the underlying non-conservative forces acting on the oscillator. The second key aspect is that the EOM is identified with only knowledge of the mass of the oscillator and the transient response. The first phase of the EDDI method constructs the dissipated energy from the kinetic energy, then identifies a mathematical model for the damping based on the dissipated energy. To achieve this, the moments in time when the displacements are zero, where the mechanical and kinetic energies are equal, are used to compute the energy dissipated by the damping of the system. The second phase begins by computing the conservative force acting on the oscillator from either a balance of the other forces in the system or through the Lagrange equation. Finally, the stiffness model is determined by solving a set of linear equations to construct a mathematical model for the conservative (elastic) force. The governing equations are discovered by incorporating both the damping and stiffness terms. The method is demonstrated by employing analytical and real measured responses of nonlinear SDOF systems with different clearances nonlinearities, which shows that the proposed approach is suitable for non-smooth mechanical systems as well as smooth systems.

## 1. Introduction

Traditionally, the governing equations that describe the dynamics of vibrating structures have been derived from the first principles and insights of the analyst, which are then validated by experimental testing [1]. Notably, the widely used experimental modal analysis [2] performs well when dealing with linear systems that have proportional or light damping and well-separated modes. However, these approaches often struggle to capture and interpret nonlinear phenomena [3]. In this sense, there is room for new methods to understand, model, and control nonlinearities in complicated systems [4], ultimately leading to improved designs, enhanced performance, and increased reliability in various engineering and scientific applications [5].

Nowadays, the availability of affordable and efficient sensors has resulted in a significant increase in measurement data, presenting the opportunity to create mathematical models using data-driven techniques [6]. One of the most effective approaches in this regard is system identification, which uses measured data (e.g., position time series) to develop and refine a model that accurately represents the dynamics of the system [7]. Nonlinear system identification (NSI) is a challenging yet crucial task in many disciplines. For example, nonlinear systems can exhibit complex behaviors, such as nonstationary, bifurcations, internal resonances, to name a few. Thus, techniques developed for linear systems may not be effective making their identification more difficult compared to linear systems [8]. Hence, NSI may require careful consideration of model structures, estimation techniques, and validation approaches.

In the field of system identification, three main categories of approaches are recognized [3]: (1) parametric methods, which assume that the underlying system can be described by a parsimonious mathematical model with a specific form, such as a differential equation. The aim is to estimate the values of the parameters that best fit the observed data [9–17] While parametric methods are powerful because they provide a clear mathematical representation of the system, and can be used for analysis and prediction, these methods require a good understanding of the governing dynamics of the system to choose an appropriate model structure [5]. (2) Non-parametric methods that do not assume a specific mathematical model of the system. Instead, they rely on the data itself to infer the behavior of the system [11, 18–24]. These methods are suitable for systems where the dynamics are complex or unknown. However, they often result in dense mathematical models, especially with noisy data. (3) Semi-parametric methods blend the elements of both parametric and non-parametric approaches, offering a balance between model structure and flexibility to account for data-driven methodologies [25–31]. These methods are particularly valuable when some insights about the data are known.

Very recently, [31] proposed a new semi-parametric method for NSI of strongly nonlinear single-degree-of-freedom (SDOF) mechanical systems called the energy-based dual-phase dynamics identification (EDDI) method. The EDDI method leverages the fact that both kinetic and mechanical energies are the same when the displacement of the oscillator is zero to identify a model for the dissipated energy (and its corresponding damping model). Once the damping model is identified, the mechanical and potential energies in the oscillator are estimated and the conservative force is computed. A system of linear equations is then constructed between a library of candidate functions and the conservative force, which is then solved to identify the parameters of the stiffness model. The method derives the governing differential equation of the nonlinear system by combining the rate of change of momentum (i.e., the inertial force), the damping model, and the stiffness models.

Mechanical systems are often subject to strongly nonlinear dynamics such as clearances, hysteresis, buckling, etc., which significantly increase their complexity and their analysis [20]. Particularly, clearance nonlinearities – referring to the gaps between components – can lead to increased noise, vibrations, and complex dynamic responses such as impacts. Several approaches have been developed to discover or model this type of nonlinearity. For instance, [32] introduced a mathematical model with periodic stiffness variation and clearance nonlinearity to describe the dynamics of the steady state forced response of a gear dynamics test rig. The study employs a generalized solution methodology based on the harmonic balance method [33], allowing for the analysis of complex non-linear phenomena commonly observed in rotating machinery. Ref. [34] identified the clearance nonlinearities, in the context of a multi-story building experiment by using the nonlinear subspace identification method [35]. Ref. [36] proposed a mathematical model for local nonlinear attachments that exhibit clearance nonlinearity, by using the characteristic nonlinear system identification method [29]. In [37], a mathematical formulation was developed for a cam-follower system that incorporates clearance nonlinearity. This formulation analyzes the transient and steady-state behaviors to understand the system's dynamics under various conditions. Ref. [38] presented a clearance-type parameter identification algorithm, and it was used to identify a single clearance nonlinearity in a cantilever beam. Based on these references, understanding these phenomena is crucial for ensuring the reliable and optimal performance of dynamic systems. Thus, the objective of this research is to perform NSI in the presence of clearance nonlinearities in SDOF systems using a modified version of the EDDI method presented in [31].

## 2. The Energy-based Dual-phase Dynamics identification (EDDI) Method

In the following, we briefly describe the EDDI method, which is shown in Fig. 1. For this method, the user needs to know the mass of the oscillator and have transient free-response data to calculate the kinetic energy of the system. When dealing with experimental measurements, such as from accelerometers, numerical integration, and filtering techniques can be applied to obtain the velocities and displacements of the oscillator [39]. In a nutshell, the method involves two phases. The first phase consists of the damping-model identification, a model whose energy dissipation corresponds to that of the measurements. The second phase consists of employing the obtained damping model and the rate of change of the momentum (e.g., for the Duffing oscillator: $\dot{p} = m\ddot{x}$) to compute the conservative force acting on the oscillator. Then, a stiffness model is determined by solving a system of linear equations formed by the generalized displacement and the conservative force. As a result, the EOM is the combination of the time derivative of the momentum, the damping model, and the stiffness model. In the following, we describe the core equations of each phase of the method.

### 2.1. Phase One: Damping Model Identification

Consider a general Duffing oscillator [40], as depicted in Fig. 1(a), whose EOM is

$$m\ddot{x} + B(x, \dot{x}, e) + K(x, e) = F(t), \qquad (1)$$

where $B(x, \dot{x})$, $K(x)$, $F(t)$, $e$, are the internal non-conservative force, the internal conservative force, an external impulsive force (or zero forcing with non-zero initial velocity), and clearance in the structure, respectively. The EDDI method was proposed as a generalized method, in a Lagrangian framework where the instantaneous mechanical is written as

$$E(t) = E(0) - \int_0^t \dot{q}B(q,\dot{q},e)d\tau, \ t \geq 0, \tag{2}$$

where the integral term represents the energy dissipated from the system. Note that Eq. 2 expresses the mechanical energy as an initial energy minus the dissipated energy and is equivalent to the standard mechanical energy given by the sum of the kinetic and potential energy.

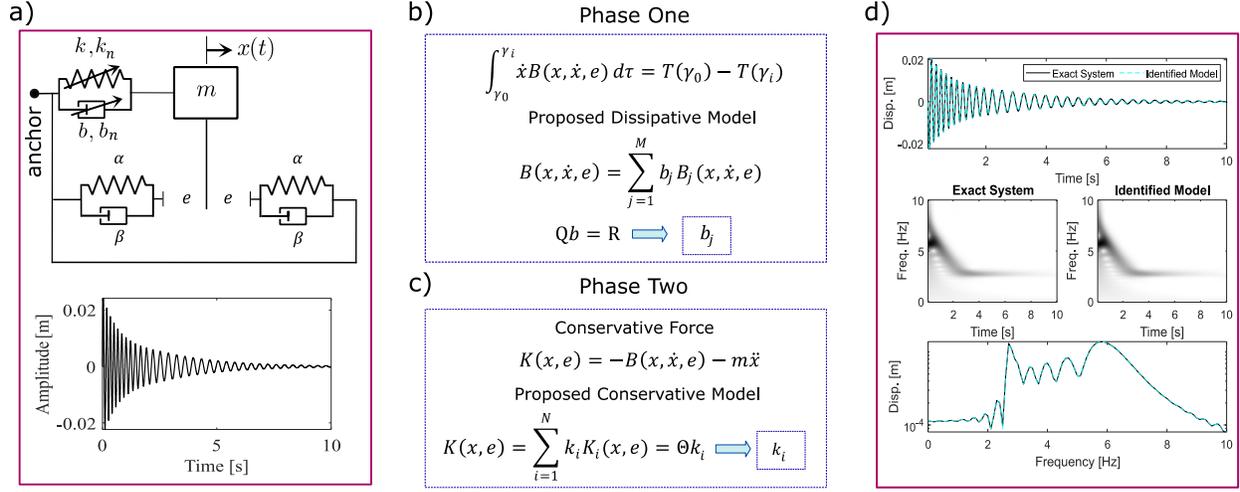

**Fig. 1.** Flowchart of the EDDI method: (a) Time series of the oscillator, (b) Phase one: dissipative-model energy matches dissipative energy in measurements and how the damping-model is obtained, (c) Phase two: the conservative force and how the stiffness model is computed, (d) validation of the resulting equation of motion.

The EDDI method takes advantage of the fact that mechanical energy is equivalent to the kinetic energy whenever the displacement is zero. We label the moments in time that satisfy this as $\gamma_i, i \in \mathbb{N}$, such that the mechanical energy is equal to the kinetic energy. Thus,

$$E(\gamma_i) = T(\gamma_i) = E(0) - \int_0^{\gamma_i} \dot{q}(\gamma_i)B(q(\gamma_i),\dot{q}(\gamma_i),e)d\tau. \tag{3}$$

For convenience, we shift time from 0 to $\gamma_0$, $E(\gamma_0) = T(\gamma_0)$, and rearranging, we get

$$\int_{\gamma_0}^{\gamma_i} \dot{q}B(q,\dot{q},e)d\tau = T(\gamma_0) - T(\gamma_i), \ i \in \mathbb{Z}^+. \tag{4}$$

From Eq. 4, $B(q,\dot{q},e)$ on the left-hand side contains the unknown model for the damping model of the oscillator, whereas the right-hand side is calculated entirely from the measurements. Note that Eq. 4 reveals a key part of the EDDI method: that the energy dissipated by the model must match the dissipated energy in the measurements. To identify the damping model, we expand the model as $B(q,\dot{q},e) = \sum_{j=1}^{M} b_j B_j(q,\dot{q},e)$, where each $B_j(q,\dot{q},e)$ is a candidate function that may be polynomial, trigonometric, or other nonlinear functions as chosen by the analyst with a scalar coefficient $b_j$ and $M$ is the number of proposed functions. It is worth mentioning that to obtain the EOM of the system, proposing a suitable library of candidate functions for the damping model is necessary. By subsituting the candidate functions in Eq. 4, we obtain a set of algebraic equations that can be written as a matrix form and truncated to $i \in [1, N]$ such that

$$\begin{bmatrix} \int_{\gamma_0}^{\gamma_1} \dot{q} B_1(q,\dot{q},e) d\tau & \cdots & \int_{\gamma_0}^{\gamma_1} \dot{q} B_M(q,\dot{q},e) d\tau \\ \vdots & \ddots & \vdots \\ \int_{\gamma_0}^{\gamma_N} \dot{q} B_1(q,\dot{q},e) d\tau & \cdots & \int_{\gamma_0}^{\gamma_N} \dot{q} B_M(q,\dot{q},e) d\tau \end{bmatrix} \begin{bmatrix} b_1 \\ \vdots \\ b_M \end{bmatrix} = \begin{bmatrix} T(\gamma_0) - T(\gamma_1) \\ \vdots \\ T(\gamma_0) - T(\gamma_N) \end{bmatrix}. \quad (5)$$

Note that Eq. 5 can be expressed in a compact form as

$$Qb = R. \quad (6)$$

Solving this equation for $b$, we obtain the parameters for the dissipation model contained in $B(q, \dot{q}, e)$.

## 2.2. Phase Two: Stiffness Model Identification

In the original EDDI method [31], the authors derived a data-driven version of the conservative force using the Lagrange equation and Lagrangian mechanics. Although successful, the approach relies on computing a numerical derivative of the Lagrangian with respect to the displacement, which results in numerical singularities when the displacement reaches a maximum or minimum. To overcome this, the original EDDI method employed a smoothing operation to remove the singularities from the computed conservative force, which although successful, slightly altered the resulting data. Here, we present another approach to calculating the conservative forces that avoids the issues that arise with computing the numerical derivatives. Instead of computing the mechanical and potential energies, we leverage the identified dissipation model, $B(q, \dot{q}, e)$, and the equation of motion similar to the restoring force method [18], such that the internal generalized conservative force is computed in terms of the other force components as

$$K(q, e) = F(t) - B(q, \dot{q}, e) - \dot{p}. \quad (7)$$

Note that $p$ can be obtained through our assumptions that the generalized coordinate and the mass of the system are known (measured). Just as the damping model is obtained, we construct the following system of equations for the conservative force using the displacement of the oscillator

$$K(q, e) = \Theta(q, e) k_i. \quad (8)$$

where $\Theta$ contains the candidate functions that will serve to identify a stiffness-model (e.g., in Cartesian coordinates: $x, x^3, (|x| - e)\text{sgn}(x)H(|x| - e)$, etc.). By solving Eq. 8, the $k_i$ coefficients that are best expressed in $\Theta(q)$ provide the internal conservative force $K(q, e)$. We refer to this form of calculating the conservative force in the EDDI method as the modified EDDI (mEDDI) method.

## 3. SI for SDOF nonlinear oscillators with a clearance nonlinearity

### 3.1. Analytical system

We demonstrate the mEDDI method in the following analytical Duffing oscillator (DO) with clearance nonlinearity with whose EOM is

$$m\ddot{x} + b\dot{x} + b_{nl}x^2\dot{x} + \beta\dot{x}H(|x| - e) + kx + k_{nl}x^3 + \alpha(|x| - e)\text{sgn}(x)H(|x| - e) = 0, \qquad (9)$$

$$x(0) = 0, \dot{x}(0) = 1,$$

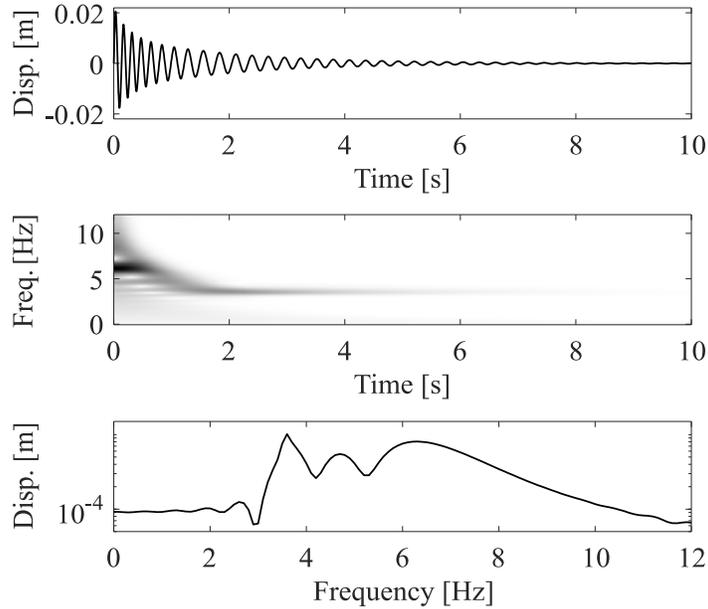

**Fig. 2.** Simulated displacement response of the DO with clearance nonlinearity.

where $m = 0.1$ kg, $b = 0.08$ Ns/m, $b_{nl} = 2000$ Ns/m³, $\beta = 0.2$ Ns/m, $k = 40$ N/m, $k_{nl} = 5000$ N/m³, $\alpha = 200$ N/m, $e = 0.005$ m and $H(\cdot)$ are the mass, linear damping, nonlinear damping, linear stiffness, nonlinear stiffness coefficients, clearance nonlinearity, the Heaviside function, respectively, and $x$ is the displacement coordinate. We use MATLAB® for all calculations throughout this paper. The simulated time series shown at the top of Fig. 2 were generated by solving Eq. 9 using the *ode45* solver with relative and absolute tolerances set to $10^{-12}$ and $10^{-16}$. Such stringent tolerances are used to ensure that the numerical integration of Eq. 9 has high accuracy given the non-smooth nature of the clearance nonlinearity. A time span of $t \in [0,10]$ s is selected, allowing the response to decay to near-zero amplitude, and the sampling rate is set to 20 kHz. We depict the continuous wavelet transform (CWT) spectrum (normalized to have a maximum amplitude of 1) [41, 42] and the Fourier spectrum in the middle and the bottom of Fig. 2, respectively. In the time-frequency representation, at a particular time and frequency the darker the shading the higher the energy content. Note that to achieve similar time-frequency representations, other techniques can also be used [43][44][45]. The CWT spectrum reveals the non-stationary nature of the response with a dark horizontal band around 6 Hz appearing for the first 0.5 seconds, which corresponds to the clearance nonlinearity.

Figure 3(a) depicts the exact linear and nonlinear damping forces within the system, indicating that the nonlinear damping is significant with the selected parameters. Along with the kinetic energy, $T(t)$, Fig. 3(b) shows the instants when the displacement is zero, $T(\gamma_i)$, as green circles. To identify the dissipative force, we propose the following model

$$B(x, \dot{x}, e) = b_1\dot{x} + b_2\dot{x}^2 + b_3\dot{x}^3 + b_4 x^2 \dot{x} + b_5 \dot{x} H(x - e) + b_6 \dot{x} H(|x| - e), \qquad (10)$$

where the terms that contain $H(\cdot)$ are necessary to capture the dynamics of the clearance nonlinearities, while the mixed displacement-velocity damping term appears when viscous dampers are excited in the transverse direction [46]. To obtain the parameter values of Eq. 10, we proceed to solve Eq. 6 by using the backslash operator (\) in MATLAB®. Table 1 provides a comparison between the calculated and the exact values. To provide comparison with other methods, we also present the values computed by the SINDy method [47] ($\lambda = 0.05$) in Table 1. The SINDy model was identified using the same candidate functions as in Eq. 10. We can see that the SINDy method performs well in identifying the damping model, but the model identified by the EDDI method is closer to the exact parameters. Note that the values for $b_2$, $b_3$ and $b_5$ by EDDI are all small enough to be neglected and would be removed if sparsity conditions were applied like

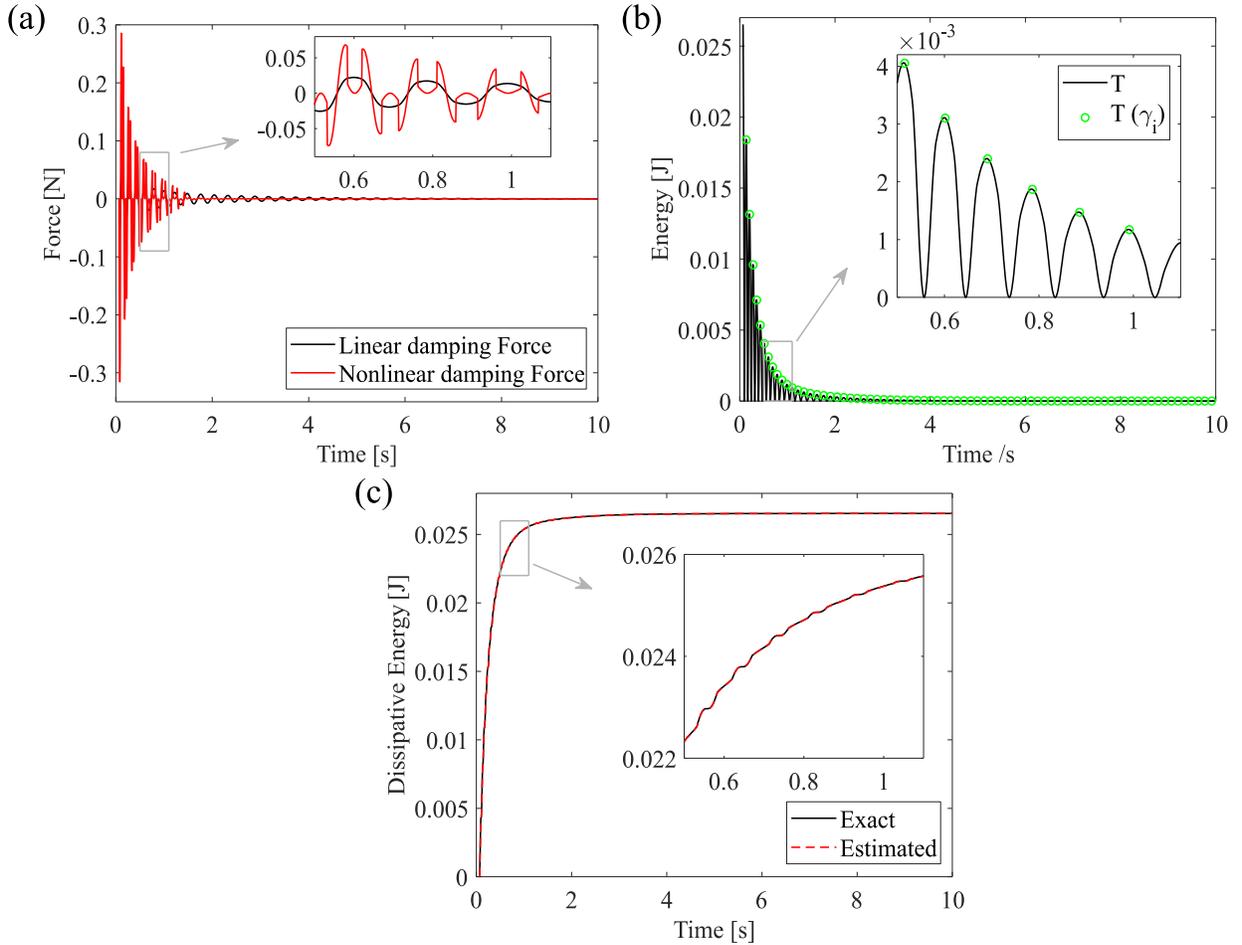

**Fig. 3.** (a) Linear and nonlinear damping forces, (b) Real kinetic and mechanical energies, (c) Comparison of the exact dissipated energy and that estimated by EDDI. (d) Comparison of the exact and estimated mechanical energies.

SINDy. Fig. 3(c) shows the exact dissipated energy from the system and that predicted by the identified damping model. The dissipated energy predicted by the damping model is obtained by integrating Eq. 10. A close agreement is observed between the two dissipated energies, which indicates that the identified damping model is accurate. Additionally, we can see that there is an initial sharp dissipation of energy that lasts about 1 second followed by a much more gradual dissipation for the rest of the response.

Table 1. Comparison of the coefficients of the dissipative model.

| Coefficient | Exact | EDDI Identified | Error | SINDy Identified | Error |
|---|---|---|---|---|---|
| $b_1 = b$ [Ns/m] | 0.08 | 0.07996 | 0.05% | 0.08016 | 0.2% |
| $b_2$ [Ns$^2$/m$^2$] | 0 | 0.00472 | – | 0 | 0 |
| $b_3$ [Ns$^3$/m$^3$] | 0 | −0.00208 | – | 0 | 0 |
| $b_4 = b_{nl}$ [Ns/m$^3$] | 2000 | 2017.1 | 0.85% | 2003.3 | 0.16% |
| $b_5$ [Ns/m] | 0 | −0.00081 | – | 0 | 0 |
| $b_6 = \beta$ [Ns/m] | 0.2 | 0.19982 | 0.09% | 0.19926 | 0.37% |

For the second phase, we know that for this system the conjugate momentum is $p = m\dot{x}$, such that we can solve for the conservative force using Eq. 7. Figure 4 shows the resulting restoring force versus the displacement computed using Eq. 7. The restoring force clearly reveals a non-smooth, clearance-type nonlinearity in the system. Thus, we choose the following stiffness model for the system

$$K(x,e) = \sum_{i=1}^{N=5} k_i x^i + k_6(x-e)H(x-e) + k_7(|x|-e)\text{sgn}(x)H(|x|-e), \qquad (11)$$

which is similar to the one proposed in [36] with the addition of the clearance nonlinearities. To identify the unknown parameters of Eq. 11, we solve Eq. 8 by using the backslash operator (\) in MATLAB®. Table 2 shows the exact and identified parameters using the EDDI, mEDDI, and SINDy methods as well as the corresponding percentage error for each parameter. We can see that the parameters $k_1 = k$, $k_3 = k_{nl}$, and $k_7 = \alpha$ are identified accurately and are similar for the original EDDI, the mEDDI, and the SINDy methods. Note that the mEDDI method provides a better estimation of the $k_1$ and $k_3$ parameters the original EDDI method, and a closer value of $k_3$ than SINDy. Plotting the identified model and the restoring force for the simulated response in the domain $x \in [-0.02, 0.02]$ m, Fig. 4 depicts the piecewise linearity behavior due to the clearance nonlinearity. At the clearance value $e$, the internal conservative force suddenly changes its behavior due to the additional stiffness that is added by the clearance nonlinearity. Furthermore, the identified stiffness model is accurate, since we can see a good agreement between the force calculated by the EDDI method and the exact force from the exact EOM.

Table 2. Comparison of the coefficients of the conservative model.

| Coefficient | Exact | original EDDI | | mEDDI | | SINDy | |
|---|---|---|---|---|---|---|---|
| | | Identified | Error | Identified | Error | Identified | Error |
| $k_1 = k$ [N/m] | 40 | 40.004 | 0.011% | 40.003 | 0.007% | 39.999 | 0.002% |
| $k_2$ [N/m$^2$] | 0 | $-8.6831$ | – | $-8.6512$ | – | 0.0153 | – |
| $k_3 = k_{nl}$ [N/m$^3$] | 5000 | 4839.4 | 3.212% | 4965.4 | 0.692% | 4957.8 | 0.845% |
| $k_4$ [N/m$^4$] | 0 | $2.58 \times 10^4$ | – | $2.54 \times 10^4$ | – | $-114.69$ | – |
| $k_5$ [N/m$^5$] | 0 | $4.79 \times 10^5$ | – | $1.67 \times 10^5$ | – | $7.28 \times 10^4$ | – |
| $k_6$ [N/m] | 0 | 0.0669 | – | 0.0652 | – | 0 | – |
| $k_7 = \alpha$ [N/m] | 200 | 199.99 | 0.004% | 199.98 | 0.010% | 200.01 | 0.004% |

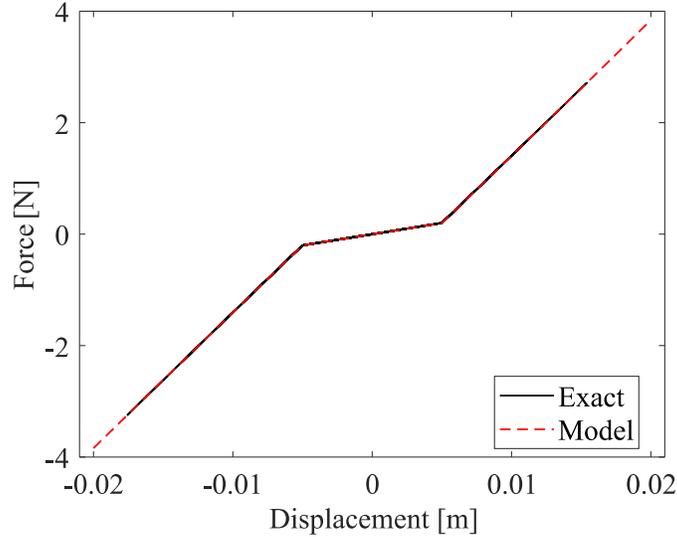

**Fig. 4.** Comparison of the internal conservative forces computed using the exact model and the model identified with the EDDI method shown as restoring force curves.

To compare the response identified system with the exact system, we combine the identified damping and stiffness models and simulate the resulting EOM using the same initial conditions. Figure 5 shows for the exact system and the identified model the displacement time series, with the corresponding CWT spectra, and Fourier spectra. We observe a strong agreement between the response of the identified model and the exact system. We anticipated a close match since this simulated signal was used for the identification process. As such, we validate the SI process by using simulated signals with different initial conditions. As shown in Fig. 6, the good agreement between the simulated signals and those obtained in this procedure validates and demonstrates the capabilities of the EDDI method to deal with clearance nonlinearities.

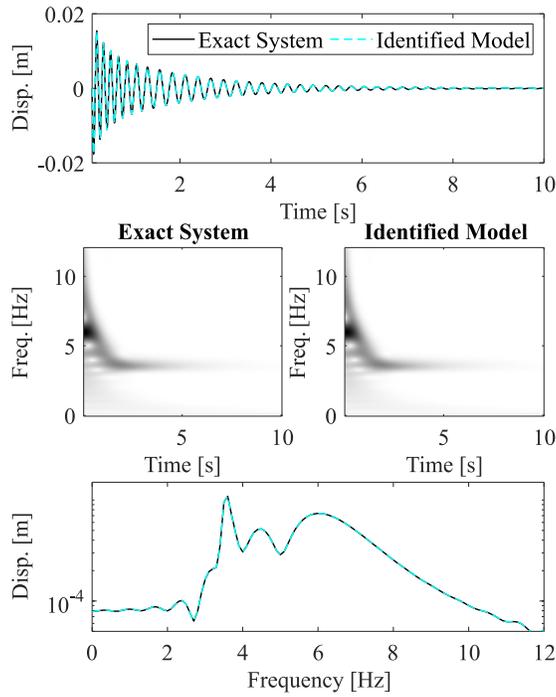

**Fig. 5.** Comparison of the displacement responses, CWT spectra, and Fourier spectra for the exact system and the identified model.

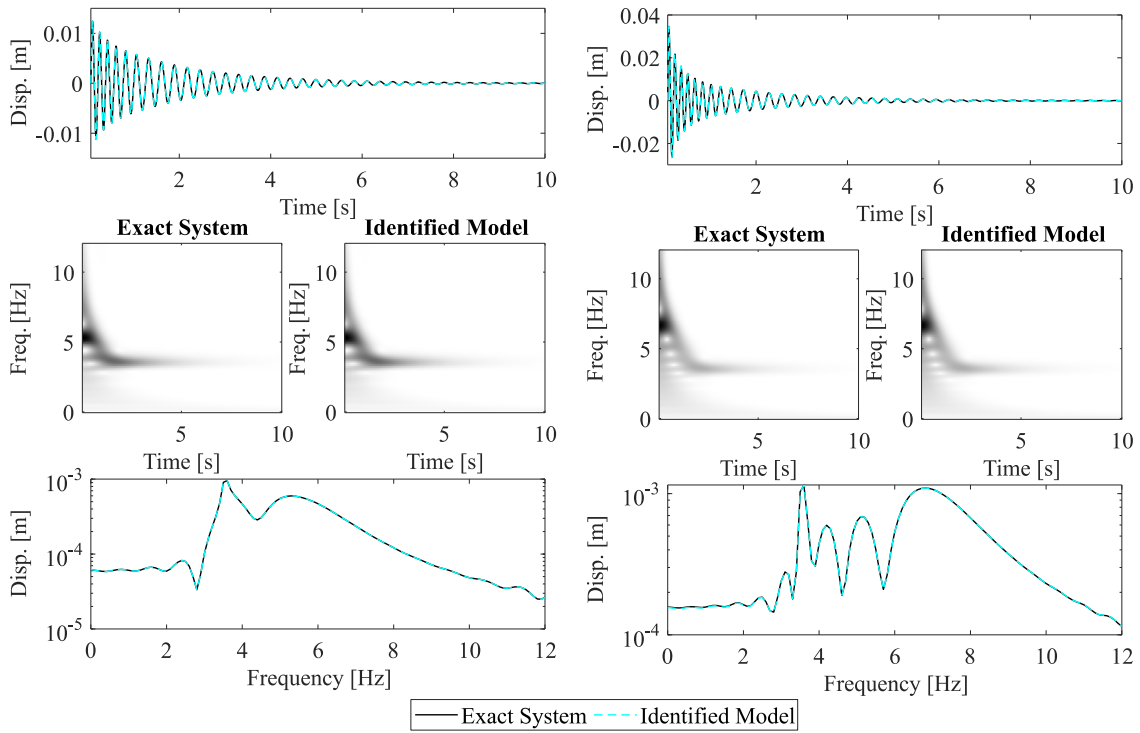

**Fig. 6.** Comparison of the simulated time series and predicted responses for initial conditions $x(0) = 0$ m, (a) $\dot{x}(0) = 0.5$ m/s and (b) $\dot{x}(0) = 2$ m/s, which were not used in the identification.

## 3.3. Experimental Demonstration

### 3.3.1 Duffing Oscillator with Clearance Nonlinearity

To demonstrate the EDDI method using a real experimental system, we considered the DO shown in Fig. 7 that uses steel flexures that introduce damping and stiffness nonlinearities. The mass of the DO is a cuboid made from steel with dimensions of 0.0305 m × 0.0305 m × 0.0127 m and a mass of 0.088 kg. The mass is suspended by (interior) steel flexures using 6-32 UNC screws that act as linear grounding springs. The flexures are anchored to a rectangular steel beam, which is bolted to an aluminum plate. The flexures have an active length of 0.0584 m, a width of 0.0102 m, and a thickness of 0.000127 m. The clearance nonlinearities are introduced by similar outer flexures, with a thickness of 0.000508, that are located at 5.35 mm from the inner flexures. This design is intentionally made symmetric, to reduce the complexity of mathematical representation. The clearance gap is realized using a vertical beam made of UHMW plastic between the two flexures at the grounding location. For simplicity, we will assume the flexures are massless. The aluminum plate is bolted with 10-32 UNC screws to L-brackets, which are grounded to the optical table. We conducted several experiments by applying impacts to the DO using an automatic modal hammer created with a stepper motor and a PCB Piezotronics modal impact hammer (model 086C03) [48][49]. The resulting motion was measured using a PCB accelerometer (model 352C03), which has a nominal sensitivity of 1 mV/(m/s$^2$). This accelerometer is robust to noise at low frequencies, which makes it a good choice for this experiment. We recorded the transient responses for 12 seconds with a sampling rate of 19200 Hz using HBM MX1601B hardware and catman data acquisition software. To get the velocities and displacements, we numerically integrated one and two times, respectively, and then high-pass filtered the outputs by using a third-order Butterworth filter with a cutoff frequency of 1.5 Hz. Note that the velocity was high-pass filtered prior to the second numerical integration, then the output was high-pass filtered again to obtain the displacements.

The external force applied to the system is shown in Fig. 8(a) and the resulting displacement response, the CWT spectrum, and the Fourier spectrum are displayed in Fig. 8(b). Due to the presence of clearance nonlinearities, the system exhibits a strongly nonlinear response with a natural frequency that begins around 6 Hz and then decays to around 4 Hz. The corresponding instantaneous kinetic energy, $T(t)$, and the kinetic energy when the displacements are zero, $T(\gamma_i)$, are provided in Fig. 8(c).

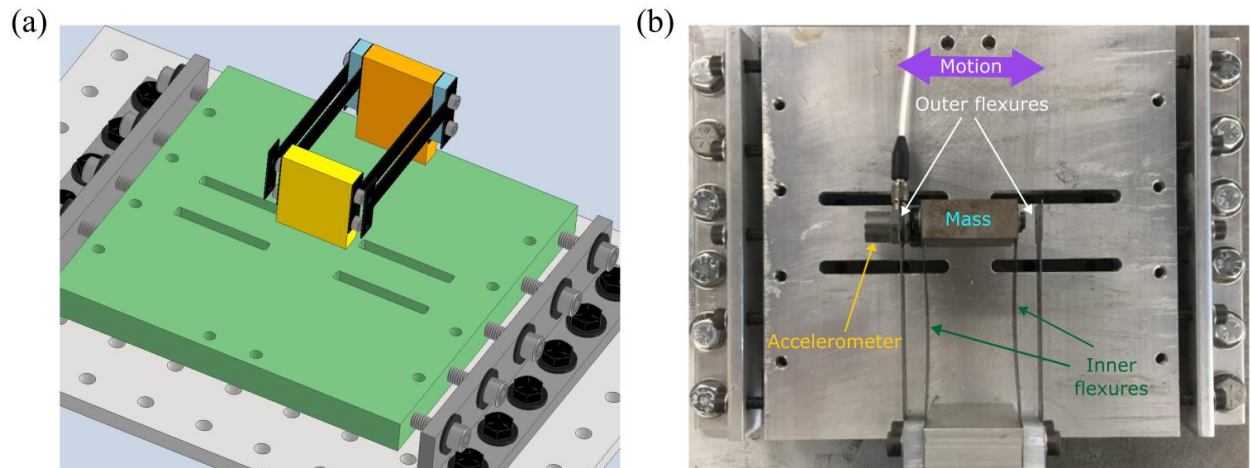

**Fig. 7.** (a) CAD model, (b) Top view of the instrumented Duffing oscillator with flexures, mass, and the direction of motion for the nonlinear oscillator.

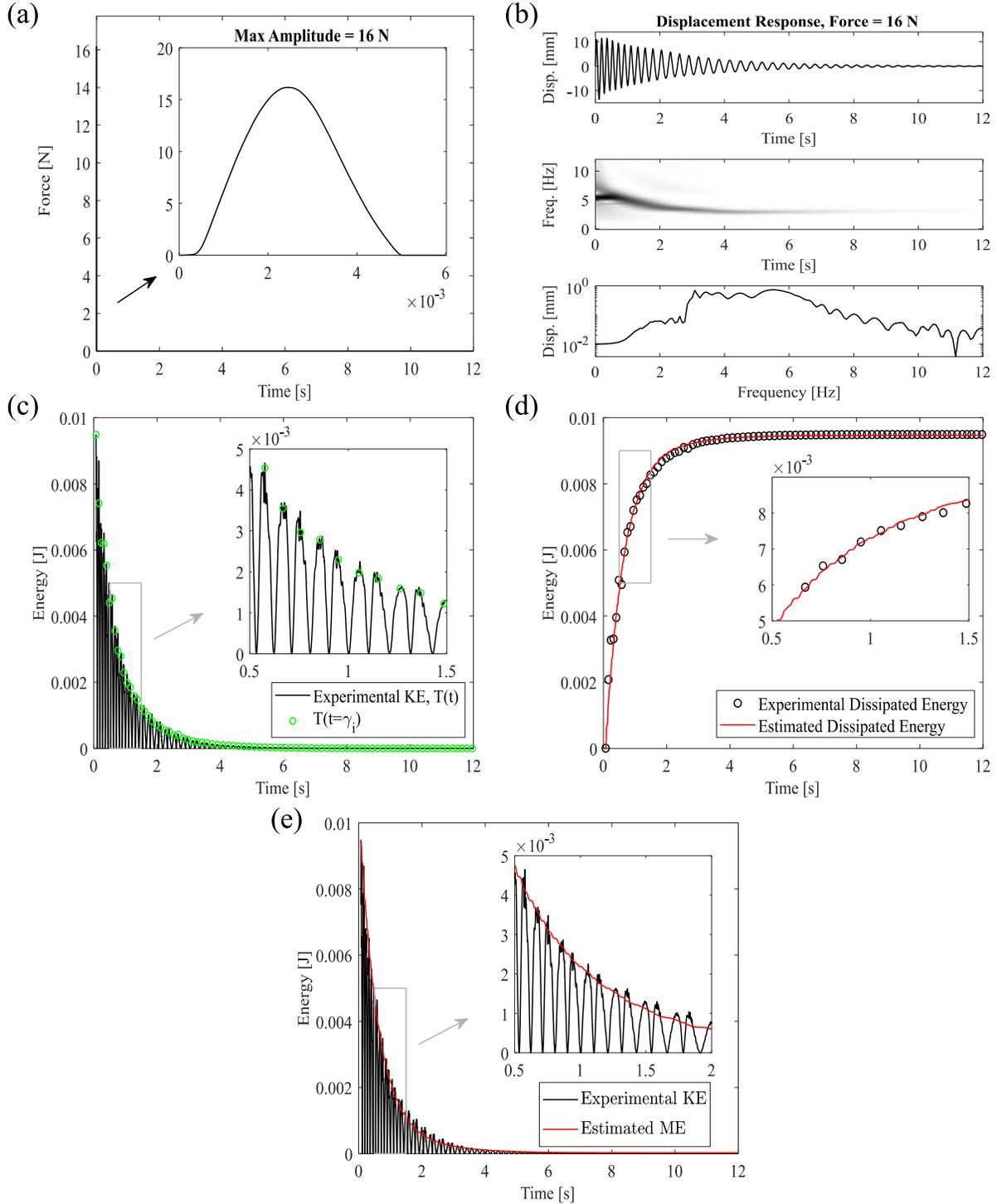

**Fig. 8.** (a) Experimentally measured excitation, (b) Displacement response of the DO-ES, (c) experimental kinetic energy and $T(\gamma_i)$, (d) Experimental and modeled dissipative energies, (e) Calculated mechanical energy.

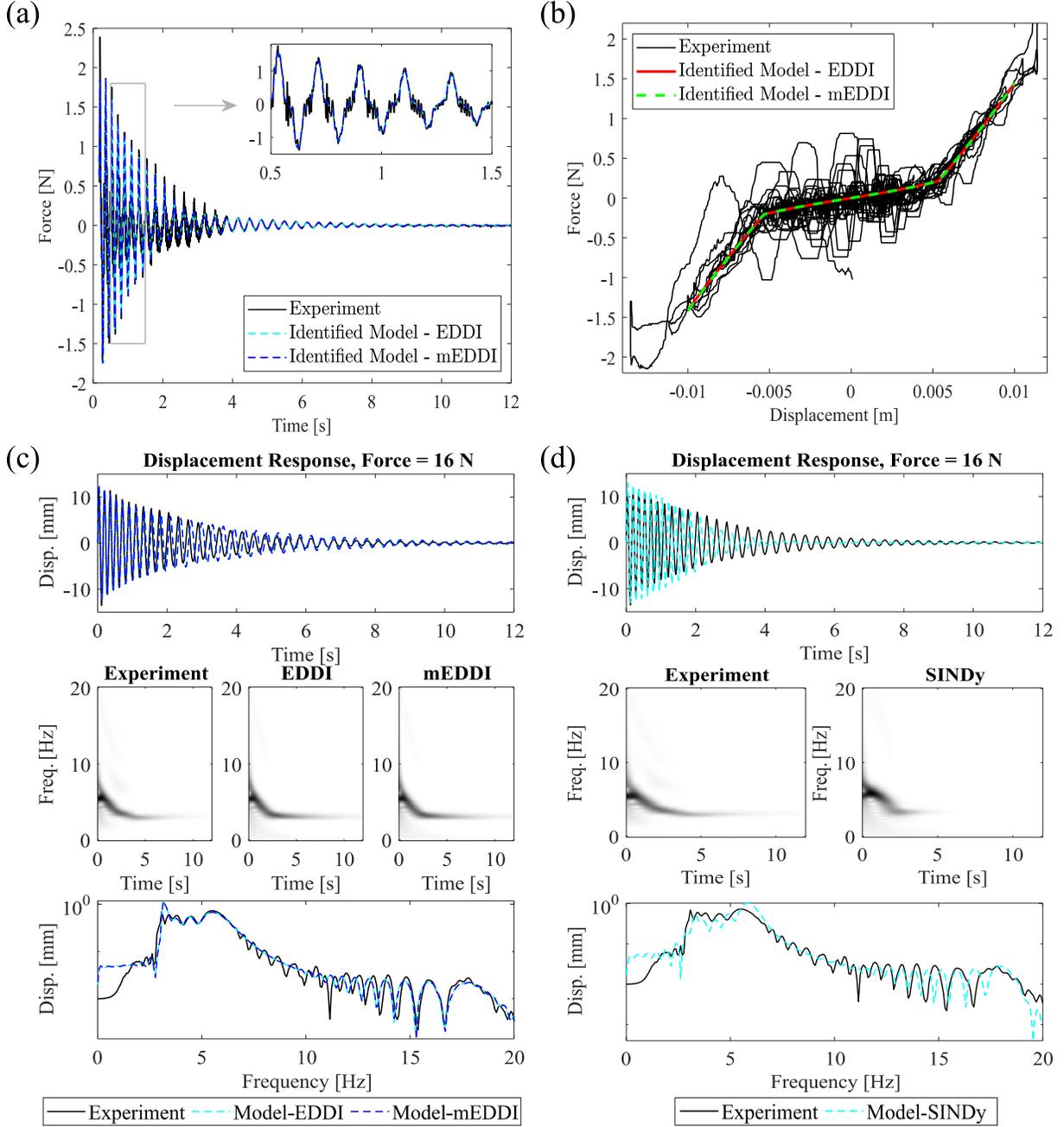

**Fig. 9.** Comparison of the internal conservative forces computed using the EDDI method (experiment) and the identified model shown as (a) time series and (b) restoring-force plot. Comparison of the measured and simulated responses for (c) the model identified using EDDI, mEDDI, and (d) the model identified by SINDy.

For the first phase, we propose the following dissipation model

$$B(x, \dot{x}, e) = b_1 \dot{x} + b_2 \dot{x} H(|x| - e). \tag{12}$$

The solution of Eq. 6 yields the following results: $b_1 = 0.056$ Ns/m, and $b_2 = 0.146$ Ns/m. Figure

8(d) shows the estimated dissipated energy. There are small discrepancies between the experimental and model dissipated energies, which arise because the model does not perfectly capture all the underlying physics as well as small effects due to noise. The mechanical energy presented in Fig. 8(e) is estimated using Eq. 2 between the initial instant, $\gamma_0$, and the final instants, $\gamma_N$. The mechanical energy shows a reasonable behavior to capture this clearance effect when the mass interacts with them.

For the second phase, we first set $p = m\dot{x}$, such that $\dot{p} = m\ddot{x}$. We then obtain the conservative force using Eq. 7 and depict it in Fig. 9(a), which displays the time series of both the experimental and model conservative forces, showing a satisfactory agreement between them. After that, we utilize the stiffness model

$$K(x,e) = \sum_{i=1}^{N=3} k_n x^i + k_4(|x| - e)\text{sgn}(x)H(|x| - e). \tag{13}$$

By solving Eq. 8, the coefficient values calculated for both the original and the mEDDI methods are similar, as presented in Table 3.

Table 3. Comparison of the coefficients of the conservative model for the real experiment

| Coefficient | original EDDI Identified | mEDDI Identified | SINDy Identified |
|---|---|---|---|
| $k_1$ [N/m] | 33.7 | 33.7 | 39.0 |
| $k_2$ [N/m$^2$] | 158.5 | 145.5 | 558.3 |
| $k_3$ [N/m$^3$] | $1.98 \times 10^5$ | $1.83 \times 10^5$ | $-1.80 \times 10^5$ |
| $\alpha$ [N/m] | 191.3 | 195.8 | 249.79 |

Applying SINDy, we get the following parameters for the dissipative model: $b_1 = 0.196$ Ns/m and $b_2 = -0.289$ Ns/m. While for the stiffness model, the calculated parameters are presented in Table 3. Fig. 9(b) depicts the restoring-force plots, and one can see the natural behavior of this system under a clearance-type, hardening phenomenon, and the EDDI method captures it. The resulting governing equation of motion is

$$m\ddot{x} + b_1\dot{x} + b_2\dot{x}H(|x| - e) + \sum_{i=1}^{N=3} k_n x^i + k_4(|x| - e)\text{sgn}(x)H(|x| - e) = F(t). \tag{14}$$

Figure 9(c) demonstrates the validation of Eq. 14 by showing a relatively good match between the measured time series and the equation's solution, and similar signals by using both EDDI and mEDDI methods. This suggests that the proposed damping and stiffness models with the identified parameters by mEDDI effectively represent the dynamics of the measured system. In Fig. 9(d), we plot the results provided by the SINDy method. As can be seen, the model does not reflect the behavior of the clearance nonlinearities correctly.

To further validate the robustness of the EDDI method, we compare simulated responses with measured responses not used in the identification. Figures 10(a) and (b) show the validation of two measurements for impacts of 9 N and 37 N, respectively. In both cases, we observe a strong agreement in the time series and the CWT spectra between the simulated responses of the model provided by the EDDI method and the measured responses. These results confirm the strength of the EDDI method to identify the dynamics of systems that have clearance nonlinearities directly from time-series measured responses.

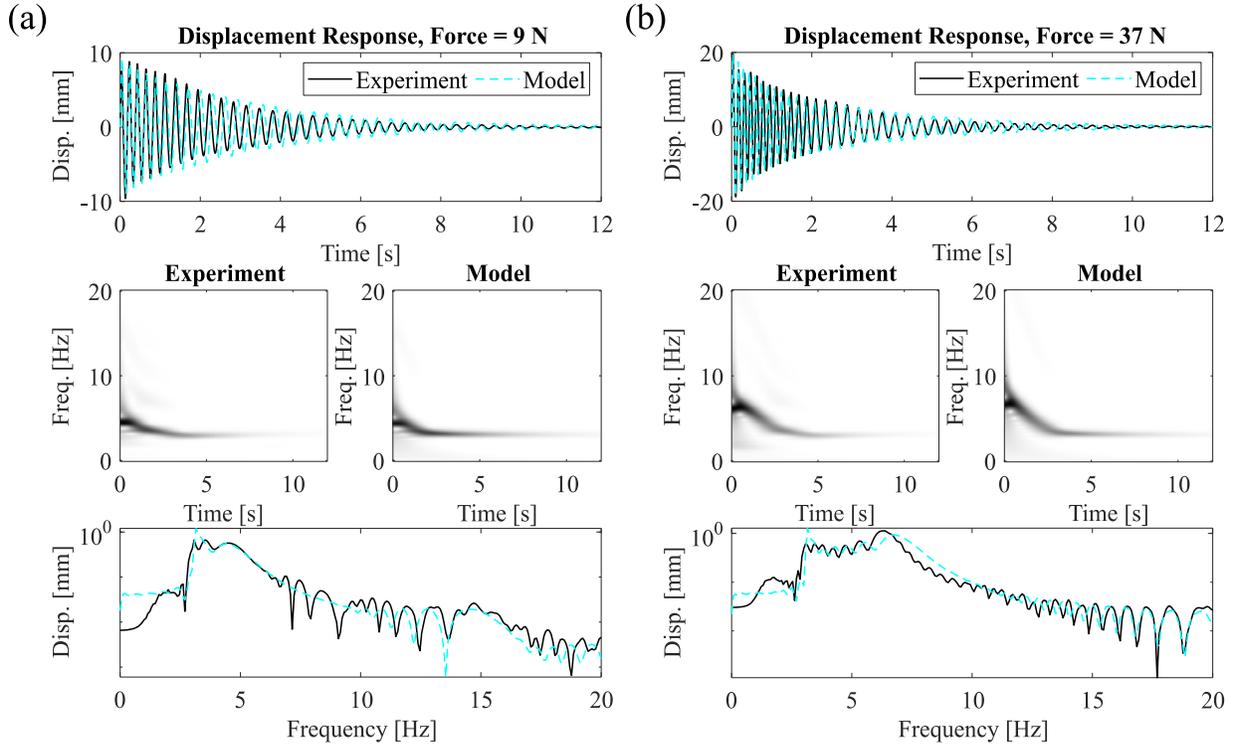

**Fig. 10.** Experimentally measured excitation and displacement responses of the DO for impact amplitudes of (a) 9 N and (b) 37 N.

## 4. Concluding Remarks

This research focused on the identification of mathematical models for SDOF systems with clearance nonlinearities in both simulated and experimentally measured responses. To achieve this aim, we used the energy-based dual-phase dynamics identification (EDDI) method, which relies only on the mass and the measured transient response. Note that no prior knowledge of the dynamics of the system is necessary in the EDDI method. The EDDI technique consists of two phases. Phase one focuses on identifying the dissipated energy and a dissipation model using the kinetic energy. In phase two, we solve for the conservative force by balancing the dissipated force, inertial force, and external force (as needed). This approach is different than the original EDDI method where the Lagrangian and a smoothing function are used to compute the conservative force. This new way of calculating the conservative forces is more convenient that the original method without changing the results substantially. Once the conservative force is obtained, a stiffness model is proposed and the parameters are identified by solving a linear system of equations between the conservative force and the library of candidate functions. Results from both simulated and experimental systems show that the EDDI method provides mathematical models that accurately represent the dynamics of SDOF systems with clearance nonlinearities. As discussed in our previous work [31], there is plenty of room to extend the EDDI technique. For instance, in the case of forced-response data, in the presence of non-smooth nonlinearities, and in multi-DOF systems, the EDDI technique would be further exploited.

**Acknowledgments**

This research was supported by the Air Force Office of Scientific Research Young Investigator Program under grant number FA9550-22-1-0295.

**Author contributions** C.L.: Writing–Original Draft, Methodology, Software, Validation, Formal analysis, Investigation, Data curation, Conceptualization. K.J.M.: Writing–Original Draft, Methodology, Software, Validation, Formal analysis, Supervision, Software, Validation, Investigation, Data curation, Project administration, Funding acquisition, Conceptualization.

**Data availability** Data sets generated during the current study are available from the corresponding author on reasonable request.

**Conflict of interest** The author declares that there is no conflict of interest or competing interest regarding the publication of this article.

326–337 (2010). https://doi.org/10.1016/j.strusafe.2010.03.00811. Frank Pai, P.: Time-frequency analysis for parametric and non-parametric identification of nonlinear dynamical systems. Mech. Syst. Signal Process. 36, 332–353 (2013). https://doi.org/10.1016/j.ymssp.2012.12.002
12. Ben Abdessalem, A., Dervilis, N., Wagg, D., Worden, K.: Model selection and parameter estimation in structural dynamics using approximate Bayesian computation. Mech. Syst. Signal Process. 99, 306–325 (2018). https://doi.org/10.1016/j.ymssp.2017.06.017
13. Ben Abdessalem, A., Dervilis, N., Wagg, D., Worden, K.: Model selection and parameter estimation of dynamical systems using a novel variant of approximate Bayesian computation. Mech. Syst. Signal Process. 122, 364–386 (2019). https://doi.org/10.1016/j.ymssp.2018.12.048
14. Scheel, M., Kleyman, G., Tatar, A., Brake, M.R.W., Peter, S., Noël, J.P., Allen, M.S., Krack, M.: Experimental assessment of polynomial nonlinear state-space and nonlinear-mode models for near-resonant vibrations. Mech. Syst. Signal Process. 143, 106796 (2020). https://doi.org/10.1016/j.ymssp.2020.106796
15. Le Guisquet, S., Amabili, M.: Identification by means of a genetic algorithm of nonlinear damping and stiffness of continuous structures subjected to large-amplitude vibrations. Part I: single-degree-of-freedom responses. Mech. Syst. Signal Process. 153, 107470 (2021). https://doi.org/10.1016/j.ymssp.2020.107470
16. Wang, C., Moore, K.J.: On nonlinear energy flows in nonlinearly coupled oscillators with equal mass. Nonlinear Dyn. 103, 343–366 (2021). https://doi.org/10.1007/s11071-020-06120-5
17. Lejarza, F., Baldea, M.: Data-driven discovery of the governing equations of dynamical systems via moving horizon optimization. Sci. Rep. 12, 1–15 (2022). https://doi.org/10.1038/s41598-022-13644-w
18. Masri, S.F., Caughey, T.K.: A nonparametric identification technique for nonlinear dynamic problems. J. Appl. Mech. Trans. ASME. 46, 433–447 (1979). https://doi.org/10.1115/1.3424568
19. Noël, J.P., Renson, L., Kerschen, G.: Complex dynamics of a nonlinear aerospace structure: Experimental identification and modal interactions. J. Sound Vib. 333, 2588–2607 (2014). https://doi.org/10.1016/j.jsv.2014.01.024
20. Vakakis, A.F., Bergman, L.A., McFarland, D.M., Lee, Y.S., Kurt, M.: Current efforts towards a non-linear system identification methodology of broad applicability. Proc. Inst. Mech. Eng. Part C J. Mech. Eng. Sci. 225, 2497–2515 (2011). https://doi.org/10.1177/0954406211417217
21. Masri, S.F., Smyth, A.W., Chassiakos, A.G., Caughey, T.K., Hunter, N.F.: Application of neural networks for detecting of changes in nonlinear systems. J. Eng. Mech. 126, 666–676 (2000)
22. Kwarta, M., Allen, M.S.: NIXO-Based identification of the dominant terms in a nonlinear equation of motion of structures with geometric nonlinearity. J. Sound Vib. 568, (2024). https://doi.org/10.1016/j.jsv.2023.117900
23. Schmidt, M., Lipson, H.: Distilling free-form natural laws from experimental data. Science (80-. ). 324, 81–85 (2009). https://doi.org/10.1126/science.1165893
24. Liu, Z., Tegmark, M.: Machine Learning Conservation Laws from Trajectories. Phys. Rev. Lett. 126, 180604 (2021). https://doi.org/10.1103/PhysRevLett.126.180604
25. Allen, M.S., Sumali, H., Epp, D.S.: Piecewise-linear restoring force surfaces for semi-
17